\documentclass[a4paper,11pt]{amsart}

\usepackage{amsmath,amssymb,amsfonts,latexsym}

\newtheorem{theorem}{Theorem}[section]

\newtheorem{corollary}[theorem]{Corollary}

\begin{document}

\title[]{A study of degenerate Bernoulli and Euler numbers via operators $\left(x^{1-\lambda}\frac{d}{dx}\right)^{n}$}

\author{Taekyun Kim}
\address{Department of Mathematics, Kwangwoon University, Seoul 139-701, Republic of Korea}
\email{tkkim@kw.ac.kr}

\author{Daesan Kim}
\address{Department of Mathematics, Sogang University, Seoul 121-742, Republic of Korea}
\email{dskim@sogang.ac.kr}



\begin{abstract}
The paper introduces novel classes of operators, $\left(x^{1-\lambda}\frac{d}{dx}\right)^{n}$ and $x^{n \lambda}\left(x^{1-\lambda}\frac{d}{dx}\right)^{n}$, and investigates their applications to the study of special numbers. We utilize these operators to derive explicit expressions for degenerate Euler numbers and type 2 degenerate Euler numbers. Furthermore, we establish a relationship between degenerate Bernoulli numbers and degenerate Euler numbers. We note here that these operators emerge naturally when we study such explicit expressions and such a relationship.
 \end{abstract}

\maketitle

\markboth{\centerline{\scriptsize A study of degenerate Bernoulli and Euler numbers via operators $\left(x^{1-\lambda}\frac{d}{dx}\right)^{n}$}}
{\centerline{\scriptsize Taekyun Kim and Dae San Kim}}


\section{Introduction}
The study of degenerate versions of certain special polynomials and numbers, which was initiated by Carlitz's works on the degenerate Bernoulli and Euler polynomials (see [4,5]), has recently seen renewed interest among mathematicians (see [3,7-9,11,13-15,21-24,26,28,29,32-34] and the references therein). They have been explored by employing various tools, including operator theory, combinatorial methods, generating functions, umbral calculus, $p$-adic analysis, probability theory, special functions, differential equations and analytic number theory. In fact, not only degenerate versions of many special polynomials and numbers have been studied but also probabilisitc extensions (see [16-10,30,37] and the references therein) and $\lambda$-analogues of them (see [12,25] and the references therein).\par
The aim of this paper is to introduce the class of operators, $\left(x^{1-\lambda}\frac{d}{dx}\right)^{n}$ (see \eqref{-2}), and to find some of their applications. These operators emerge naturally in connection with explicit determination of the degenerate Euler numbers $\mathcal{E}_{n,\lambda}$ (see \eqref{8}) and finding a relationship between the degenerate Bernoulli numbers $\beta_{n,\lambda}$ (see \eqref{7}) and the degenerate Euler numbers. In addition, the class of operators, $x^{n\lambda}\left(x^{1-\lambda}\frac{d}{dx}\right)^{n}$ (see \eqref{-1}), are used in explicit determination of the type 2 degenerate Euler numbers $E_{n,\lambda}$ (see \eqref{9}).  \par
In more detail, the outline of this paper is as follows. We introduce and study one class of operators given by
\begin{equation}
\left(x^{1-\lambda}\frac{d}{dx}\right)^{n}=\sum_{l=0}^{n}{n \brace l}_{\lambda}x^{l-n\lambda}\left(\frac{d}{dx}\right)^{l}, \label{-2}
\end{equation}
and another class of operators given by
\begin{equation}
x^{n\lambda}\left(x^{1-\lambda}\frac{d}{dx}\right)^{n}=\sum_{l=0}^{n}{n \brace l}_{\lambda}x^{l}\left(\frac{d}{dx}\right)^{l}=\left(x\frac{d}{dx}\right)_{n,\lambda}. \label{-1}
\end{equation}
Application of \eqref{-1} to $f(x)=e^{x}$ yields the Dobinski-like formula in Theorem 2.3. We get the identity in Theorem 2.4 by applying \eqref{-1} to $\frac{x^{a}}{1-x^{b}}$, for any real numbers $a,b$ with $b \ne 0$. This makes us introduce the polynomials $S_{n,\lambda}(x;a,b)$ (see \eqref{18}), which are cooked up from the degenerate Frobenius polynomials $F_{n,\lambda}(x)$ (see \eqref{5}). In Theorem 2.5, we obtain an explicit expression of the degenerate Euler numbers involving the degenerate Stirling numbers of the second kind ${n \brace k}_{\lambda}$ (see \eqref{1}). This is done by evaluating $\left(\frac{d}{dz}\right)^{n}\big(\frac{1}{1+e_{\lambda}(z)}\big)$ in two different ways at $z=0$, where $e_{\lambda}(z)$ are the degenerate exponentials in \eqref{0}. In one way of evaluation, \eqref{-2} is used by observing that $\left(\frac{d}{dz}\right)^{n}\big(\frac{1}{1+e_{\lambda}(z)}\big)=\left(u^{1-\lambda}\frac{d}{du}\right)^{n}\left(\frac{1}{1+u}\right)$, with the substitution $u=e_{\lambda}(z)$ (see \eqref{19}). In Theorem 2.7, we derive a relationship between the degenerate Bernoulli numbers and the degenerate Euler numbers. This again involves the evaluation of $\left(\frac{d}{dz}\right)^{n}\big(\frac{1}{1+e_{\lambda}(z)}\big)$ at $z=0$, which makes use of \eqref{-2}, as we mentioned in the above. Let $g_{n,\lambda}(u)=u^{n\lambda}\left(u^{1-\lambda}\frac{d}{du}\right)^{n}\left(\frac{u}{1+u^{2}}\right)$. Then we show in Theorem 2.9 that $g_{n,\lambda}(u)$ is equal to $\sum_{k=0}^{n}a(n,k|\lambda)u^{2k+1}(u^2+1)^{-k-1}$, with explicitly determined $a(n,k|\lambda)$ as a finite sum. Finally, we determine the type 2 degenerate Euler numbers $E_{n,\lambda}$ by evaluating $e_{\lambda}^{n\lambda}(x)\frac{d^{n}}{dx^{n}}\mathrm{sech}_{\lambda}(x)$ at $x=0$ in two different ways, where $\mathrm{sech}_{\lambda}(x)$ is the degenerate hyperbolic secant function defined in \eqref{10}. Here again one way of the evaluation uses the observation that $e_{\lambda}^{n\lambda}(x)\frac{d^{n}}{dx^{n}}\mathrm{sech}_{\lambda}(x)=2g_{n,\lambda}(u)$, with the substitution $u=e_{\lambda}(z)$ (see \eqref{36}). Thus the value of $e_{\lambda}^{n\lambda}(x)\frac{d^{n}}{dx^{n}}\mathrm{sech}_{\lambda}(x)$ at $x=0$
is $2g_{n,\lambda}(1)=\sum_{k=0}^{n}a(n,k|\lambda)2^{-k}$. So, the explicit determination of $E_{n,\lambda}$ makes use of the operators in \eqref{-1}. As general references, the reader may refer to [1,2,6,10,27,31,35,36]. For the rest of this section, we recall the facts that are needed throughout this paper. \par

\vspace{0.1in}

For any nonzero $\lambda\in \mathbb{R}$, the degenerate exponentials are defined by
\begin{equation}
e_{\lambda}^{x}(t)=\sum_{n=0}^{\infty}(x)_{n,\lambda}\frac{t^{n}}{n!}, \ e_{\lambda}(t)=e_{\lambda}^{1}(t), \ (\mathrm{see}\ [11,15,23]), \label{0}
\end{equation}
where the degenerate falling factorials $(x)_{n,\lambda}$ are given by
\begin{equation*}
(x)_{0,\lambda}=1, (x)_{n,\lambda}=x(x-\lambda)(x-2\lambda)\cdots(x-(n-1)\lambda), \ (n \geq 1).
\end{equation*}
Note that $\lim_{\lambda \to 0}e_{\lambda}^{x}(x)=e^{xt}$. \par
The degenerate Stirling numbers of the second kind are defined by
\begin{equation}
(x)_{n,\lambda}=\sum_{k=0}^{n}{n \brace k}_{\lambda}(x)_{k}, (n \geq 0), \ (\mathrm{see}\ [11]), \label{1}
\end{equation}
where the falling factorials are given by
\begin{equation*}
(x)_{0}=1, \ (x)_{n}=x(x-1)\cdots(x-n+1), \ (n \geq 1).
\end{equation*}
Note that $\lim_{\lambda \to 0}{n \brace k}_{\lambda}={n \brace k}$,
where ${n \brace k}$ are the Stirling numbers of the second kind defined by
\begin{equation}
x^{n}=\sum_{k=0}^{n}{n \brace k}(x)_{k}, \ (n \geq 0),\ (\mathrm{see}\ [1,6,35,36]).\nonumber \\ \nonumber
\end{equation}
From \eqref{1}, we note that
\begin{equation}
\frac{1}{k!}(e_{\lambda}(t)-1)^{k}=\sum_{n=k}^{\infty}{n \brace k}_{\lambda}\frac{t^{n}}{n!},\  (k \geq 0), \ (\mathrm{see}\ [11,16,21-26]). \label{2}
\end{equation} \par
We recall from [22,25] that
\begin{equation}
\left(x\frac{d}{dx}\right)_{n,\lambda}=\sum_{k=0}^{n}{n \brace k}_{\lambda}x^{k}\left(\frac{d}{dx}\right)^{k}, \ (n \geq 0). \label{3}
\end{equation}
The degenerate Fubini polynomials are given by
\begin{equation}
\frac{1}{1-x(e_{\lambda}(t)-1)}=\sum_{k=0}^{\infty}F_{k,\lambda}(x)\frac{t^{k}}{k!}, \ (\mathrm{see}\ [37]). \label{4}
\end{equation}
By \eqref{2} and \eqref{4}, we get
\begin{equation}
F_{n,\lambda}(x)=\sum_{k=0}^{n}{n \brace k}_{\lambda}k!x^{k},\ (n \geq 0), \ (\mathrm{see}\ [37]). \label{5}
\end{equation}
From \eqref{3}, we note that
\begin{align}
\sum_{k=0}^{\infty}(k)_{n,\lambda}x^{k}&=\left(x\frac{d}{dx}\right)_{n,\lambda}\left(\frac{1}{1-x}\right)=\frac{1}{1-x}\sum_{k=0}^{n}{n \brace k}_{\lambda}k!\left(\frac{x}{1-x}\right)^{k} \label{6} \\
&\ =\frac{1}{1-x}F_{n,\lambda}\left(\frac{x}{1-x}\right),\ (n \geq 0), \  (\mathrm{see}\ [37]). \nonumber
\end{align}
By taking $\lambda \rightarrow 0$ of \eqref{6}, we get
\begin{equation*}
\sum_{k=0}^{\infty}k^{m}x^{k}=\left(x\frac{d}{dx}\right)^{m}\left(\frac{1}{1-x}\right)=\frac{1}{1-x}F_{m}\left(\frac{x}{1-x}\right), \ (m \geq 0), \ (\mathrm{see}\ [27]),
\end{equation*}
where $F_{m}(x)$ are the Fubini polynomials defined by
\begin{equation*}
F_{m}(x)=\sum_{k=0}^{m}{m \brace k}k!x^{k},\ (m \geq 0).
\end{equation*} \par
In [4,5], Carlitz introduced the degenerate Bernoulli numbers given by
\begin{equation}
\frac{t}{e_{\lambda}(t)-1}=\sum_{n=0}^{\infty}\beta_{n,\lambda}\frac{t^{n}}{n!}. \label{7}
\end{equation}
Note that $\lim_{\lambda \to 0}\beta_{n,\lambda}=B_{n}, \ (n \geq 0)$, where $B_{n}$ are the Bernoulli numbers defined by
\begin{equation}
\frac{t}{e^{t}-1}=\sum_{n=0}^{\infty}B_{n}\frac{t^{n}}{n!}, (\mathrm{see}\ [2,6,35,36]). \nonumber \\ \nonumber
\end{equation}
The degenerate Euler numbers are given by
\begin{equation}
\frac{2}{e_{\lambda}(t)+1}=\sum_{n=0}^{\infty}\mathcal{E}_{n,\lambda}\frac{t^{n}}{n!},\ (\mathrm{see}\ [4,5]). \label{8}
\end{equation}
Note that $\lim_{\lambda \to 0}\mathcal{E}_{n,\lambda}=E_{n}, \ (n \geq 0)$, where $E_{n}$ are the Euler numbers defined by
\begin{equation}
\frac{2}{e^{t}+1}=\sum_{n=0}^{\infty}E_{n}\frac{t^{n}}{n!}. \nonumber \\ \nonumber
\end{equation}
The type 2 degenerate Euler numbers are given by
\begin{equation}
\frac{2}{e_{\lambda}(t)+e_{\lambda}^{-1}(t)}=\sum_{n=0}^{\infty}E_{n,\lambda}\frac{t^{n}}{n!},\ (\mathrm{see}\ [13,14,29]). \label{9}
\end{equation}
Note that $\lim_{\lambda \to 0}E_{n,\lambda}=E_{n}^{*}, \ (n \geq 0)$, are the type 2 Euler numbers defined by
\begin{equation}
\frac{2}{e^{t}+e^{-t}}=\sum_{n=0}^{\infty}E_{n}^{*}\frac{t^{n}}{n!}, \ (\mathrm{see}\ [1,6,35]).\nonumber \\ \nonumber
\end{equation} \par
For any $\lambda \in \mathbb{R}$, we consider the degenerate hyperbolic cosine and hyperbolic sine functions which are respectively given by
\begin{equation*}
\cosh_{\lambda}(t)=\frac{e_{\lambda}(t)+e_{\lambda}^{-1}(t)}{2},\quad \sinh_{\lambda}(t)=\frac{e_{\lambda}(t)-e_{\lambda}^{-1}(t)}{2}.
\end{equation*}
In addition, we define the degenerate hyperbolic secant function as
\begin{equation}
\mathrm{sech}_{\lambda}(t)=\frac{1}{\cosh_{\lambda}(t)}=\frac{2}{e_{\lambda}(t)+e_{\lambda}^{-1}(t)}. \label{10}
\end{equation} \par
The degenerate Bell polynomials are given by
\begin{equation}
\phi_{n,\lambda}(x)=\sum_{k=0}^{n}{n \brace k}_{\lambda}x^{k}, \ (n \geq 0), \ (\mathrm{see}\ [10,11,36]). \label{11}
\end{equation}
Note that $\lim_{\lambda \to 0}\phi_{n,\lambda}(x)=\sum_{k=0}^{n}{n \brace k}x^{k}=\phi_{n}(x)$,
where $\phi_{n}(x)$ are the Bell polynomials.

\section{A study of degenerate Bernoulli and Euler numbers via operators $\left(x^{1-\lambda}\frac{d}{dx}\right)^{n}$}

For any integers $n,k $ with $n,k \ge 0$, and any $\lambda \in \mathbb{R}$, we observe that
\begin{align}
\left(x^{1-\lambda}\frac{d}{dx}\right)^{n}x^{k}&=k(k-\lambda)\cdots(k-(n-1)\lambda)x^{k-n\lambda} \label{12}\\
&\ =(k)_{n,\lambda}x^{k-n\lambda}=\sum_{l=0}^{n}{n \brace l}_{\lambda}x^{l-n\lambda}\left(\frac{d}{dx}\right)^{l}x^{k}. \nonumber
\end{align}
For any formal power series $f(x)=\sum_{k=0}^{\infty}a_{k}x^{k}$, from \eqref{12} we have
\begin{equation}
    \left(x^{1-\lambda}\frac{d}{dx}\right)^{n}f(x)=\sum_{l=0}^{n}{n \brace l}_{\lambda}x^{l-n\lambda}\left(\frac{d}{dx}\right)^{l}f(x). \label{13}
\end{equation}
By \eqref{13}, we get
\begin{equation}
    x^{n\lambda}\left(x^{1-\lambda}\frac{d}{dx}\right)^{n}f(x)=\sum_{l=0}^{n}{n \brace l}_{\lambda}x^{l}\left(\frac{d}{dx}\right)^{l}f(x). \label{14}
\end{equation}
Therefore, from \eqref{3} and \eqref{14} we obtain the following theorem.
\begin{theorem}
For any integer $n \ge 0$, we have
\begin{equation*}
x^{n\lambda}\left(x^{1-\lambda}\frac{d}{dx}\right)^{n}f(x)=\sum_{l=0}^{n}{n \brace l}_{\lambda}x^{l}\left(\frac{d}{dx}\right)^{l}f(x)=\left(x\frac{d}{dx}\right)_{n,\lambda}f(x).
\end{equation*}
\end{theorem}
From \eqref{13}, we note that
\begin{align}
\left(x^{1-\lambda}\frac{d}{dx}\right)^{n+1}f(x)&=\left(x^{1-\lambda}\frac{d}{dx}\right)\left(x^{1-\lambda}\frac{d}{dx}\right)^{n}f(x) \label{15} \\
&=x^{1-\lambda}\frac{d}{dx}\sum_{l=0}^{n}{n \brace l}_{\lambda}x^{l-n\lambda}\left(\frac{d}{dx}\right)^{l}f(x) \nonumber \\
&\ =\sum_{l=0}^{n}{n \brace l}_{\lambda}(l-n\lambda)x^{l-(n+1)\lambda}\left(\frac{d}{dx}\right)^{l}f(x)\nonumber \\
&\quad\quad+\sum_{l=0}^{n}{n \brace l}_{\lambda}x^{l+1-(n+1)\lambda}\left(\frac{d}{dx}\right)^{l+1}f(x) \nonumber \\
&\ =\sum_{l=0}^{n}{n \brace l}_{\lambda}(l-n\lambda)x^{l-(n+1)\lambda}\left(\frac{d}{dx}\right)^{l}f(x)\nonumber \\
&\quad\quad+\sum_{l=1}^{n+1}{n \brace l-1}_{\lambda}x^{l-(n+1)\lambda}\left(\frac{d}{dx}\right)^{l}f(x). \nonumber
\end{align}
Thus, by \eqref{15}, we get
\begin{equation}
 x^{(n+1)\lambda}\left(x^{1-\lambda}\frac{d}{dx}\right)^{n+1}f(x)=\sum_{l=0}^{n+1}\left({n \brace l}_{\lambda}(l-n\lambda)+{n \brace l-1}_{\lambda}\right)x^{l}\left(\frac{d}{dx}\right)^{l}f(x). \label{16}
\end{equation}
Therefore, by \eqref{16} and Theorem 2.1, we obtain the following theorem.
\begin{theorem}
For any integer $n \geq 0$, we have
\begin{displaymath}
\sum_{l=0}^{n+1}{n+1 \brace l}_{\lambda}x^{l}\left(\frac{d}{dx}\right)^{l}f(x)=\sum_{l=0}^{n+1}\left({n \brace l}_{\lambda}(l-n\lambda)+{n \brace l-1}_{\lambda}\right)x^{l}\left(\frac{d}{dx}\right)^{l}f(x).
\end{displaymath}
\end{theorem}
\noindent In particular, for $1 \leq l \leq n$, we have
\begin{equation}
{n+1 \brace l}_{\lambda}={n \brace l}_{\lambda}(l-n\lambda)+{n \brace l-1}_{\lambda}.\nonumber \\ \nonumber
\end{equation}

\vspace{0.1in}

Taking $f(x)=e^{x}$ in Theorem 2.1 and from \eqref{11}, we get the following Dobinski-like formula.
\begin{theorem}
For any integer $n \geq 0$, we have
\begin{displaymath}
 e^{x}\phi_{n,\lambda}(x)=e^{x}\sum_{k=0}^{n}{n \brace k}_{\lambda}x^{k}=\sum_{j=0}^{\infty}\frac{(j)_{n,\lambda}}{j!}x^{j}.
\end{displaymath}
\end{theorem}
Taking $f(x)=\frac{1}{1-x}$, by \eqref{6} and Theorem 2.1, we have
\begin{equation*}
x^{n\lambda}\left(x^{1-\lambda}\frac{d}{dx} \right)^{n}\left(\frac{1}{1-x}\right)=\frac{1}{1-x}F_{n,\lambda}\left(\frac{x}{1-x}\right)=\sum_{k=0}^{\infty}(k)_{n,\lambda}x^{k}.
\end{equation*}
More generally, for $a,b \in \mathbb{R}$ with $b \neq 0$, and any integer $n \geq 0$, we have
\begin{equation}
x^{n\lambda}\left(x^{1-\lambda}\frac{d}{dx}\right)^{n}\left(\frac{x^{a}}{1-x^{b}}\right)=x^{n\lambda}\left(x^{1-\lambda}\frac{d}{dx}\right)^{n}\sum_{k=0}^{\infty}x^{a+bk}. \label{17}
\end{equation}
{\it{We define, for $a,b \in \mathbb{R}$ with $b \neq 0$ and any integer $n \geq 0$,  the polynomial $S_{n,\lambda}(x:a,b)$ by
\begin{equation}
S_{n,\lambda}(x;a,b)=b^{n}\sum_{k=0}^{n}\binom{n}{k}\left(\frac{a}{b}\right)_{n-k,\frac{\lambda}{b}}F_{k,\frac{\lambda}{b}}(x), \label{18}
\end{equation}
where $F_{k,\frac{\lambda}{b}}(x)$ is the degenerate Fubini polynomial given in \eqref{5}.}} \\
\vspace{0.1in}
Note that $S_{n,\lambda}(x;0,1)=F_{n,\lambda}(x), \ (n \geq 0)$.
Then, by \eqref{17} and \eqref{18}, we get
\begin{align}
&x^{n\lambda}\left(x^{1-\lambda}\frac{d}{dx}\right)^{n}\left(\frac{x^{a}}{1-x^{b}}\right)=x^{n\lambda}\left(x^{1-\lambda}\frac{d}{dx}\right)^{n}\sum_{k=0}^{\infty}x^{a+bk} \label{18-1} \\
&\ =\sum_{k=0}^{\infty}(a+bk)_{n,\lambda}x^{a+bk}=x^{a}b^{n}\sum_{k=0}^{\infty}\left(\frac{a}{b}+k\right)_{n,\frac{\lambda}{b}}x^{bk} \nonumber \\
&\ =x^{a}b^{n}\sum_{k=0}^{\infty}\sum_{l=0}^{n}\binom{n}{l}\left(\frac{a}{b}\right)_{n-l,\frac{\lambda}{b}}(k)_{l,\frac{\lambda}{b}}x^{bk} \nonumber \\
&=x^{a}b^{n}\sum_{l=0}^{n}\binom{n}{l}\left(\frac{a}{b}\right)_{n-l,\frac{\lambda}{b}}\sum_{k=0}^{\infty}(k)_{l,\frac{\lambda}{b}}x^{bk} \nonumber \\
&\ =x^{a}b^{n}\sum_{l=0}^{n}\binom{n}{l}\left(\frac{a}{b}\right)_{n-l,\frac{\lambda}{b}}\frac{1}{1-x^{b}}F_{l,\frac{\lambda}{b}}\left(\frac{x^{b}}{1-x^{b}}\right) \nonumber \\
&=\frac{x^{a}}{1-x^{b}}b^{n}\sum_{l=0}^{n}\binom{n}{l}\left(\frac{a}{b}\right)_{n-l,\frac{\lambda}{b}}F_{l,\frac{\lambda}{b}}\left(\frac{x^{b}}{1-x^{b}}\right) \nonumber \\
&\ =\frac{x^{a}}{1-x^{b}}S_{n,\lambda}\left(\frac{x^{b}}{1-x^{b}};a,b\right), \ (n \geq 0). \nonumber
\end{align}
Therefore, by \eqref{18-1}, we obtain the following theorem.
\begin{theorem}
For any integer $n \geq 0$ and $a,b \in \mathbb{R}$ with $b \neq 0$, we have
\begin{displaymath}
x^{n\lambda}\left(x^{1-\lambda}\frac{d}{dx}\right)^{n}\left(\frac{x^{a}}{1-x^{b}}\right)=\frac{x^{a}}{1-x^{b}}S_{n,\lambda}\left(\frac{x^{b}}{1-x^{b}};a,b\right),
\end{displaymath}
where $S_{n,\lambda}(x;a,b)$ is the polynomial given in \eqref{18}.
\end{theorem}

Let $u=e_{\lambda}(z)$. Then, by \eqref{13}, we have
\begin{align}
    \left(\frac{d}{dz}\right)^{n}\left(\frac{1}{1+e_{\lambda}(z)}\right)&=\left(\frac{d}{du}\frac{du}{dz}\right)^{n}\left(\frac{1}{1+u}\right)=\left(u^{1-\lambda}\frac{d}{du}\right)^{n}\left(\frac{1}{1+u}\right) \label{19}  \\
    &\ =\sum_{k=0}^{n}{n \brace k}_{\lambda}u^{k-n\lambda}\left(\frac{d}{du}\right)^{k}\left(\frac{1}{1+u}\right)  \nonumber \\
    &\ =\sum_{k=0}^{n}{n \brace k}_{\lambda}u^{k-n\lambda}(-1)^{k}k!(1+u)^{-k-1} \nonumber \\
    &\ =\sum_{k=0}^{n}{n \brace k}_{\lambda}e_{\lambda}^{k-n\lambda}(z)(-1)^{k}k!(1+e_{\lambda}(z))^{-k-1}. \nonumber
\end{align}
Thus, by \eqref{19}, we get
\begin{equation}
\left(\frac{d}{dz}\right)^{n}\left(\frac{1}{1+e_{\lambda}(z)}\right)\bigg\vert_{z=0}=\sum_{k=0}^{n}{n \brace k}_{\lambda}(-1)^{k}k!2^{-k-1}. \label{20}
\end{equation}
On the other hand, by \eqref{8}, we get
\begin{align}
\left(\frac{d}{dz}\right)^{n}\left(\frac{1}{1+e_{\lambda}(z)}\right)\bigg\vert_{z=0}&=\frac{1}{2}\left(\frac{d}{dz}\right)^{n}\left(\frac{2}{1+e_{\lambda}(z)}\right)\bigg\vert_{z=0} \label{21} \\
&\ =\frac{1}{2}\mathcal{E}_{n,\lambda}, \ (n \geq 0). \nonumber
\end{align}
Therefore, by \eqref{20} and \eqref{21}, we obtain the following theorem.
\begin{theorem}
For any integer $n \geq 0$, we have
\begin{displaymath}
\mathcal{E}_{n,\lambda}=\sum_{k=0}^{n}{n \brace k}_{\lambda}(-1)^{k}k!2^{-k}.
\end{displaymath}
\end{theorem}
For $n \in \mathbb{N}$, we have
\begin{align}
    &\left(\frac{d}{dz}\right)^{n}\big(zf(z)\big)=\left(\frac{d}{dz}\right)^{n-1}\left(\frac{d}{dz}\big(zf(z)\big)\right)=\left(\frac{d}{dz}\right)^{n-1}\left(f(z)+z\frac{d}{dz}f(z)\right) \label{22} \\
    &\ =\left(\frac{d}{dz}\right)^{n-1}f(z)+\left(\frac{d}{dz}\right)^{n-2}\frac{d}{dz}\left(z\frac{d}{dz}f(z)\right) \nonumber \\
    &\ =\left(\frac{d}{dz}\right)^{n-1}f(z)+\left(\frac{d}{dz}\right)^{n-2}\left(\frac{d}{dz}f(z)+z\left(\frac{d}{dz}\right)^{2}f(z)\right) \nonumber \\
    &\ =2\left(\frac{d}{dz}\right)^{n-1}f(z)+\left(\frac{d}{dz}\right)^{n-2}\left(z\left(\frac{d}{dz}\right)^{2}f(z)\right) \nonumber
\end{align}
\begin{align*}
    &\ =2\left(\frac{d}{dz}\right)^{n-1}f(z)+\left(\frac{d}{dz}\right)^{n-3}\frac{d}{dz}\left(z\left(\frac{d}{dz}\right)^{2}f(z)\right) \nonumber \\
    &\ =2\left(\frac{d}{dz}\right)^{n-1}f(z)+\left(\frac{d}{dz}\right)^{n-3}\left(\left(\frac{d}{dz}\right)^{2}f(z)+z\left(\frac{d}{dz}\right)^{3}f(z)\right)  \nonumber \\
    &\ =3\left(\frac{d}{dz}\right)^{n-1}f(z)+\left(\frac{d}{dz}\right)^{n-3}\left(z\left(\frac{d}{dz}\right)^{3}f(z)\right) \nonumber \\
    &\ =\cdots \nonumber \\
    &\ =n\left(\frac{d}{dz}\right)^{n-1}f(z)+z\left(\frac{d}{dz}\right)^{n}f(z).  \nonumber
\end{align*}
Therefore, by \eqref{22}, we obtain the following theorem.
\begin{theorem}
For any integer $n \ge 1$, we have
\begin{displaymath}
\left(\frac{d}{dz}\right)^{n}\big(zf(z)\big)=n\left(\frac{d}{dz}\right)^{n-1}f(z)+z\left(\frac{d}{dz}\right)^{n}f(z).
\end{displaymath}
\end{theorem}
From \eqref{7}, we note that
\begin{equation}
    \frac{2z}{e_{\lambda}^{2}(z)-1}=\frac{2z}{e_{\frac{\lambda}{2}}(2z)-1}=\sum_{k=0}^{\infty}2^{k}\beta_{k,\frac{\lambda}{2}}\frac{z^{k}}{k!}. \label{23}
\end{equation}
Thus, by \eqref{23}, we get
\begin{align}
\left(\frac{d}{dz}\right)^{n}\left(\frac{2z}{e_{\lambda}^{2}(z)-1}\right)\bigg\vert_{z=0}&=\left(\frac{d}{dz}\right)^{n}\left(\frac{2z}{e_{\frac{\lambda}{2}}(2z)-1}\right)\bigg\vert_{z=0} \label{24}\\
&\ =2^{n}\beta_{n,\frac{\lambda}{2}}, \quad(n \geq 0). \nonumber
\end{align}
We observe that
\begin{equation}
\frac{2z}{e_{\frac{\lambda}{2}}(2z)-1}= \frac{2z}{e_{\lambda}^{2}(z)-1}=\frac{z}{e_{\lambda}(z)-1}-\frac{z}{e_{\lambda}(z)+1}. \label{25}
\end{equation}
From \eqref{25} and Theorem 2.6, we note that
\begin{align}
&\left(\frac{d}{dz}\right)^{n}\left(\frac{2z}{e_{\frac{\lambda}{2}}(2z)-1}-\frac{z}{e_{\lambda}(z)-1}\right)=-\left(\frac{d}{dz}\right)^{n}\left(\frac{z}{e_{\lambda}(z)+1}\right) \label{26} \\
&\ =-\left(z\left(\frac{d}{dz}\right)^{n}\left(\frac{1}{e_{\lambda}(z)+1}\right)+n\left(\frac{d}{dz}\right)^{n-1}\left(\frac{1}{e_{\lambda}(z)+1}\right)\right).  \nonumber
\end{align}
Evaluating the identity \eqref{26} at $z=0$, and using \eqref{20}, \eqref{21} and \eqref{24}, we have
\begin{align}
2^{n}\beta_{n,\frac{\lambda}{2}}-\beta_{n,\lambda}&=-n\left(\frac{d}{dz}\right)^{n-1}\left(\frac{1}{e_{\lambda}(z)+1}\right)\bigg\vert_{z=0} \label{27}\\
&\ =n\sum_{k=0}^{n-1}{n-1 \brace k}_{\lambda}(-1)^{k-1}k!2^{-k-1} \nonumber \\
&=-\frac{n}{2}\mathcal{E}_{n-1,\lambda}, \ (n \in \mathbb{N}). \nonumber
\end{align}
Therefore, by \eqref{27}, we obtain the following theorem.
\begin{theorem}
For any integer $n \ge 1$, we have
\begin{displaymath}
2^{n}\beta_{n,\frac{\lambda}{2}}-\beta_{n,\lambda}=n\sum_{k=0}^{n-1}{n-1 \brace k}_{\lambda}(-1)^{k-1}k!2^{-k-1}=-\frac{n}{2}\mathcal{E}_{n-1,\lambda}.
\end{displaymath}
\end{theorem}
By taking $\lambda \rightarrow 0$ in the identities of Theorem 2.7, we obtain
\begin{corollary}
For any integer $n \ge 1$, we have
\begin{equation*}
B_{n}=-\frac{n}{2(2^{n}-1)}E_{n-1}=\frac{n}{2^{n}-1}\sum_{k=0}^{n-1}{n-1 \brace k}(-1)^{k-1}k!2^{-k-1}.
\end{equation*}
\end{corollary}
Let
\begin{equation}
g_{n,\lambda}(u)=u^{n\lambda}\left(u^{1-\lambda}\frac{d}{du}\right)^{n}\left(\frac{u}{1+u^{2}}\right), \ (n \geq 0). \label{28}
\end{equation}
Then, by \eqref{14}, we get
\begin{align}
g_{n,\lambda}(u)=u^{n\lambda}\left(u^{1-\lambda}\frac{d}{du}\right)^{n}\left(\frac{u}{1+u^{2}}\right)&=u^{n\lambda}\left(u^{1-\lambda}\frac{d}{du}\right)^{n}\sum_{j=0}^{\infty}(-1)^{j}u^{2j+1} \label{29} \\
&\ =\sum_{j=0}^{\infty}(-1)^{j}(2j+1)_{n,\lambda}u^{2j+1}. \nonumber
\end{align}
Assume that
\begin{equation}
g_{n,\lambda}(u)=\sum_{k=0}^{n}a(n,k\vert\lambda)u^{2k+1}(u^{2}+1)^{-k-1}. \label{30}
\end{equation}
Then, by binomial expansion, we get
\begin{align}
g_{n,\lambda}(u)&=\sum_{k=0}^{n}a(n,k\vert\lambda)u^{2k+1}(u^{2}+1)^{-k-1} \label{31} \\
&\ =\sum_{k=0}^{n}a(n,k\vert\lambda)u^{2k+1}\sum_{i=0}^{\infty}\binom{k+i}{k}(-1)^{i}u^{2i} \nonumber \\
&\ =\sum_{j=0}^{\infty}\left(\sum_{k=0}^{j}\binom{j}{k}a(n,k\vert\lambda)(-1)^{k-j}\right)u^{2j+1}. \nonumber
\end{align}
By \eqref{29} and \eqref{31}, we get
\begin{equation}
(2j+1)_{n,\lambda}=\sum_{k=0}^{j}\binom{j}{k}a(n,k\vert\lambda)(-1)^{k}. \label{32}
\end{equation}
From \eqref{32} and by binomial inversion (see [35], (7.12), p. 83), we have
\begin{equation}
a(n,j\vert\lambda)=\sum_{k=0}^{j}\binom{j}{k}(2k+1)_{n,\lambda}(-1)^{k}. \label{33}
\end{equation}
Therefore, by \eqref{30} and \eqref{33}, we obtain the following theorem.
\begin{theorem}
For any integer $n \ge 0$, we let
\begin{equation*}
g_{n,\lambda}(u)=u^{n\lambda}\left(u^{1-\lambda}\frac{d}{du}\right)^{n}\left(\frac{u}{1+u^{2}}\right).
\end{equation*}
Then we have
\begin{equation*}
g_{n,\lambda}(u)=\sum_{k=0}^{n}a(n,k\vert\lambda)u^{2k+1}(u^{2}+1)^{-k-1},
\end{equation*}
\end{theorem}
\noindent with
\begin{equation}
a(n,k\vert\lambda)=\sum_{j=0}^{k}(-1)^{j}\binom{k}{j}(2j+1)_{n,\lambda}.\nonumber \\ \nonumber
\end{equation} \par
From \eqref{9} and \eqref{10}, we note that
\begin{align}
e_{\lambda}^{n\lambda}(x)\frac{d^{n}}{dx^{n}}\mathrm{sech}_{\lambda}(x)&=e_{\lambda}^{n\lambda}(x)\left(\frac{d}{dx}\right)^{n}\left(\frac{2}{e_{\lambda}(x)+e_{\lambda}^{-1}(x)}\right) \label{34} \\
&\ =\sum_{l=0}^{\infty}(n\lambda)_{l,\lambda}\frac{x^{l}}{l!}\sum_{k=0}^{\infty}E_{k+n,\lambda}\frac{x^{k}}{k!}\nonumber \\
&=\sum_{m=0}^{\infty}\sum_{k=0}^{m}\binom{m}{k}E_{k+n,\lambda}(n\lambda)_{m+k,\lambda}\frac{x^{m}}{m!}.  \nonumber
\end{align}
Thus, by \eqref{34}, we get
\begin{equation}
e_{\lambda}^{n\lambda}(x)\frac{d^{n}}{dx^{n}}(\mathrm{sech}_{\lambda}(x))\Big\vert_{x=0}=E_{n,\lambda}, \ (n \geq 0). \label{35}
\end{equation}
Let $e_{\lambda}(x)=u$. Then, by \eqref{28}, we get
\begin{align}
e_{\lambda}^{n\lambda}(x)\frac{d^{n}}{dx^{n}}\mathrm{sech}_{\lambda}(x)&=e_{\lambda}^{n\lambda}(x)\left(\frac{d}{dx}\right)^{n}\left(\frac{2}{e_{\lambda}(x)+e_{\lambda}^{-1}(x)}\right) \label{36} \\
&\ =e_{\lambda}^{n\lambda}(x)\left(\frac{d}{dx}\right)^{n}\left(\frac{2e_{\lambda}(x)}{e_{\lambda}^{2}(x)+1}\right) \nonumber \\
&=2u^{n\lambda}\left(u^{1-\lambda}\frac{d}{du}\right)^{n}\left(\frac{u}{1+u^{2}}\right) \nonumber \\
&\ =2g_{n,\lambda}(u)=2g_{n,\lambda}(e_{\lambda}(x)). \nonumber
\end{align}
Taking $x=0$ in \eqref{36}, we have
\begin{align}
e_{\lambda}^{n\lambda}(x)\frac{d^{n}}{dx^{n}}\mathrm{sech}_{\lambda}(x)\Big\vert_{x=0}&=2g_{n,\lambda}(1) \label{37} \\
&\ =2\sum_{k=0}^{n}a(n,k\vert\lambda)2^{-k-1}  \nonumber \\
&\ =\sum_{k=0}^{n}a(n,k\vert\lambda)2^{-k}. \nonumber
\end{align}
Therefore, by \eqref{35}, \eqref{37} and Theorem 2.9, we obtain the following theorem.
\begin{theorem}
For any integer $n \geq 0$, we have
\begin{displaymath}
E_{n,\lambda}=\sum_{k=0}^{n}2^{-k}\sum_{j=0}^{k}(-1)^{j}\binom{k}{j}(2j+1)_{n,\lambda}.
\end{displaymath}
\end{theorem}
\section{Conclusion}
In this paper, we used the operators $\left(x^{1-\lambda}\frac{d}{dx}\right)^{n}$ in connections with explicit determination of the degenerate Euler numbers $\mathcal{E}_{n,\lambda}$ and finding a relationship between the degenerate Bernoulli numbers $\beta_{n,\lambda}$ and the degenerate Euler numbers. In addition, the operators $x^{n\lambda}\left(x^{1-\lambda}\frac{d}{dx}\right)^{n}$ are utilized in deriving an explicit expression for the type 2 degenerate Euler numbers $E_{n,\lambda}$. \par
Among other things, we showed the following:
\begin{align*}
&\mathcal{E}_{n,\lambda}=\sum_{k=0}^{n}{n \brace k}_{\lambda}(-1)^{k}k!2^{-k}, \\
&2^{n}\beta_{n,\frac{\lambda}{2}}-\beta_{n,\lambda}=n\sum_{k=0}^{n-1}{n-1 \brace k}_{\lambda}(-1)^{k-1}k!2^{-k-1}=-\frac{n}{2}\mathcal{E}_{n-1,\lambda},\\
&E_{n,\lambda}=\sum_{k=0}^{n}2^{-k}\sum_{j=0}^{k}(-1)^{j}\binom{k}{j}(2j+1)_{n,\lambda},
\end{align*}
where ${n \brace k}_{n,\lambda}$ are the degenerate Stirling numbers of the second kind.
As we noted earlier, the aforementioned operators arise naturally when we make change of variables in the relevant computations. \par
We would like to continue to investigate probabilistic extensions and $\lambda$-analogues of special polynomials and numbers as well as degenerate versions of them.

\end{document}